\documentclass{amsart}
\newtheorem{Theorem}{Theorem}[section]
\newtheorem{Proposition}[Theorem]{Proposition}
\newtheorem{Definition}[Theorem]{Definition}

\newtheorem{Lemma}[Theorem]{Lemma}
\newtheorem{Corollary}[Theorem]{Corollary}
\theoremstyle{remark}
\newtheorem{Example}[Theorem]{Example}
\def\il{\int\limits_}

\def\eps{\varepsilon}

\def\ovr{\overline}

\def\Om{\Omega}
\def\al{\alpha}

\def\Dl{\Delta}
\def\dl{\delta}

\def\bd{\partial}
\def\lm{\lambda}

\def\sm{\setminus}
\def\sbs{\subset}

\def\Sbs{\subset\subset}

\def\wtl{\widetilde}

\def\Arg{\operatorname{Arg}}

\def\im{\mathbf {Im\,}}
\def\be{\begin{enumerate}}
\def\ee{\end{enumerate}}
\def\bT{\begin{Theorem}}
\def\eT{\end{Theorem}}
\def\bP{\begin{Proposition}}
\def\eP{\end{Proposition}}
\def\bD{\begin{Definition}}
\def\eD{\end{Definition}}
\def\bE{\begin{Example}}
\def\eE{\end{Example}}
\def\bL{\begin{Lemma}}
\def\eL{\end{Lemma}}
\def\bC{\begin{Corollary}}
\def\eC{\end{Corollary}}
\def\beq{\begin{equation}}
\def\eeq{\end{equation}}

\def\aB{\mathbb B}
\def\aD{\mathbb D}
\def\aT{\mathbb T}
\def\aC{\mathbb C}

\def\aR{\mathbb R}

\def\aC{\mathbb C}

\begin{document}
\title{ (Pluri)potential compactifications}
\author{ Evgeny A. Poletsky}
\begin{abstract}  Using pluricomplex Green functions we introduce a compactification of a complex manifold $M$ invariant with respect to biholomorphisms similar to the Martin compactification in the potential theory.
\par For this we show the existence of a norming volume form $V$ on $M$ such that all negative plurisubharmonic functions on $M$ are in $L^1(M,V)$. Moreover, the set of such functions with the norm not exceeding 1 is compact. Identifying a point $w\in M$ with the normalized pluricomplex Green function with pole at $w$ we get an imbedding of $M$ into a compact set and the closure of $M$ in this set is the pluripotential compactification.
\thanks{The author was partially supported by a grant from Simons Foundation.}
\end{abstract} \keywords{Plurisubharmonic functions,
pluripotential theory, Martin boundary} \subjclass[2010]{ Primary: 32J05; secondary:
31C35, 32U15}
\address{ Department of Mathematics,  215 Carnegie Hall,
Syracuse University,  Syracuse, NY 13244}
\maketitle
\section{Introduction}
\par In this paper we construct biholomorphically invariant compactifications of complex manifolds. For domains in the complex plane there is the Carath\'eodory  compactification that is invariant with respect to biholomorphisms. It is constructed using prime ends. There are papers that used this notion for higher dimensions but it, seemingly, did not lead to invariant compactifications.
\par Our construction is similar to the Martin compactification but instead of Green functions that are not biholomorphically invariant we use their analog on complex manifolds, namely, the pluricomplex Green functions.
\par The classical Martin's approach is to consider the normalized Green functions $\wtl G_D(x,y)=G_D(x,y)/G_D(x_0,y)$ on a domain $D\sbs\aR^n$, where $x_0$ is a fixed point in $D$, and then define the Martin boundary as the set of all sequences $\wtl G_D(x,y_j)$ that converge in $L^1_{loc}$. Due to Harnack's inequalities the choice of the point $x_0$ is non-essential. The limit is a harmonic function on $D$ that is called the Martin kernel.
\par On a complex manifold $M$ this normalization does not work because pluricomplex Green functions $g_M(z,w)$ are only maximal, i. e., $(dd^c_zg_M(z,w))^n=0$ outside of $w$. There are no Harnack's inequalities and the existence of a subsequence converging in $L^1_{loc}$ is not guaranteed.
\par To circumvent this obstacle we show the existence of a norming volume form $V$ on $M$ or $D$ such that all negative (pluri)subharmonic functions on $M$ or $D$ are in $L^1(M,V)$ or $L^1(D,V)$ respectively. Moreover, the set of such functions with the norm not exceeding 1 is compact. Identifying a point $w\in M$ with the normalized pluricomplex Green function with pole at $w$ we get an imbedding of $M$ into a compact set and the closure of $M$ in this set is the pluripotential compactification. The same approach in the real case produces the Martin compactification.
\par Unfortunately, we were not able to prove that the limits of pluricomplex Green functions are maximal. It is known due to Lelong \cite{L} that almost any plurisubharmonic function is the limit in $L^1_{loc}$ of maximal functions. Thus the general theory is not applicable. At the last section we compute pluripotential compactification for a ball, smooth strongly convex domains and a bidisk. In all cases the limits are maximal and are scalar multiples of pluriharmonic Poisson kernels computed in \cite{D} and \cite{BPT}. In the two first cases the pluripotential boundary coincides with the Euclidean boundary while in the case of a bidisk it is the product of a circle and a 2-sphere.
\par We are grateful to the referee whose corrections and comments significantly improved the exposition.
\section{Green functions}
\par We denote by $G_D(x,y)$ the {\it negative} Green function on a domain $D$ in $\mathbb R^n$. It is known that the Green function is symmetric
and continuous in $G_D\times G_D$ and if $\bd D$ is $C^2$-smooth, then
(see \cite[24.1]{V})  $G_D(x,y)$ is continuous on $D\times\ovr D$ and
at every point $y\in\bd D$ there is the derivative
\begin{equation}\label{e:dgf} P_D(x,y)=\frac{\bd G_D(x,y)}{\bd n(y)},\end{equation}
along the outward normal vector $n(y)$ to $\bd D$ at $y$.
The function $P(x,y)$ is harmonic in $x$, positive on $D$ and is
called the Poisson kernel of $D$.
\par By \cite[Theorem 6.18]{H} every negative subharmonic function $u$ on $D$
can be represented as
\[u(x)=\int\limits_D G_D(x,y)\Delta u(y)+H(x),\]
where $\Dl u$ is the Riesz mass of $u$ and $H(x)$ is the least harmonic majorant of $u$. It was proved in
\cite{PK} that if $H(x)$ is a negative harmonic function on a domain $D$ with $C^2$-boundary, then
\[H(x)=-\int\limits_{\bd D}P_D(x,y)\,d\mu(y),\]
where $\mu$ is a Borel measure on $\bd D$.
\par Combining the two previous equations we get the Poisson--Jensen
formula for a negative subharmonic function $u$ on $D$:
\begin{equation}\label{e:pjb} u(x)=
\int\limits_D G_D(x,y)\Delta u(y)- \int\limits_{\bd
D}P_D(x,y)\,d\mu(y).\end{equation}
\par Let $\dl(x)$ be the distance from $x$ to $\bd D$. The following
theorem was proved in \cite{KL} and \cite[Eqns. (5) and (7)]{S} (see also \cite{W} and \cite{Z}).
\bT\label{T:kls} If $D\sbs\aR^n$ is a domain with $C^2$-boundary, then there is
a constant $A>0$ depending only on $D$ such that the Green function
satisfies the inequality
\begin{equation}\label{e:bn1}
G_D(x,y)\ge-\frac12\ln\left(1+A\frac{\dl(x)\dl(y)}{|x-y|^2}\right),\quad
n=2,
\end{equation}
and
\begin{equation}\label{e:bn2}
G_D(x,y)\ge-A\frac{\dl(x)\dl(y)}{|x-y|^n},\quad n>2.
\end{equation}\eT
\par It follows from this inequalities (see \cite{PK} and \cite{S})
that there is a constant $B>0$ depending only on $D$ such that
\begin{equation}\label{e:bn3}
P_D(x,y)\le\frac{B\dl(x)}{|x-y|^n}
\end{equation} for all $n\ge 2$.
\section{Norming functions and Martin compactification}
\par If $\phi>0$ is a continuous function on a domain $D\sbs\aR^n$, then we denote by $L^1(D,\phi)$ the space of all Lebesgue measurable functions $u$ on $D$ such that
\[\|u\|_\phi=\il D|u|\phi\,dx<\infty.\]
Let us call a positive continuous function $\phi(z)$ on $D$
{\it norming } if for each compact set $F\Sbs D$ there is a
positive constant $C(F)$ and for every increasing sequence of
subdomains $D_j\Sbs D$ with $\cup D_j=D$ there is a
sequence of numbers $\eps_j>0$ converging to zero such that for every
negative subharmonic function $u$ on $D$:
\be
\item $\|u\|_\phi<\infty$;
\item $u(x)\le -C(F)\|u\|_\phi$ for every point $x\in F$;
\item
$\|u\|_{L^1(D\sm D_j,\phi)}\le\eps_j\|u\|_\phi$.
\ee
\par The first important feature of norming functions $\phi$ on $D$ is the integrability of all negative subharmonic functions. Hence, the cone $SH^-(D)$ of all non-positive subharmonic functions on $D$  lies in the Banach space $L^1(D,\phi)$.
\par As the following lemma shows this embedding of $SH^-(D)$ into $L^1(D,\phi)$ practically does not depend on the choice of $\phi$.
\bL\label{L:en} Norming functions determine equivalent norms on
the cone of negative subharmonic functions.
\eL
\begin{proof} Let $\phi_1$ and $\phi_2$ be norming functions on
$D$. We take a compact set $F\sbs D$ of positive Lebesgue measure
and find a constant $c>0$ such that
$u(x)\le -c\|u\|_{\phi_1}$ for all $x\in F$ and all negative subharmonic functions $u$ on $D$. Then
\[-\|u\|_{\phi_2}\le\il Fu\,\phi_2\,dx\le -c\|u\|_{\phi_1}\il
F\,\phi_2\,dx.\] This shows that
there is a constant $a>0$ such that $\|u\|_{\phi_2}\ge
a\|u\|_{\phi_1}$ on $SH^-(D)$.
\end{proof}
\par Another important feature of norming functions $\phi$ is
given by the following theorem.
\bT\label{T:ct} The set $SH^-_1(D,\phi)=\{u\in SH^-(D):\,\|u\|_\phi\le1\}$ is
compact in $L^1(D,\phi)$.\eT
\begin{proof} Let  $v_j$ be a sequence of subharmonic functions in $SH^-_1(D)$. By \cite[Theorem 4.1.9]{Ho} there is a subsequence $v_{j_k}$ converging  in $L^1_{loc}(D)$ to a subharmonic function $v$ on $D$. Moreover,
\[\limsup_{k\to\infty}v_{j_k}(z)\le v(z)\]
and the left and the right side of this inequality are equal a.e.
\par By the third  property of norming functions
this subsequence converges to $v$ in $L^1(D,\phi)$ and by Fatou's lemma $\|v\|_{\phi}\le1$.
\end{proof}
\par The following lemma provides estimates of integrals of
subharmonic functions on compact sets.
\bL\label{L:ei} Let $\phi$ be a norming function on a domain
$D\sbs{\mathbb R}^n$ and let $F$ be a compact set in $D$ with $a=\|\chi_F\|_{\phi}>0$. Then there is a positive constant $c$,
depending only on $F$, such that for every negative subharmonic
function $u$ on $D$ we have
\[ca\|u\|_\phi\le-\il Fu\phi\,dx.\]\eL
\begin{proof} By the second property of norming functions
\[\int\limits_Fu\phi\,dx\le -C(F)\|u\|_\phi\il F\phi\,dx=-C(F)\|\chi_F\|_{\phi}\|u\|_\phi.\]
\end{proof}
\par Let $D\sbs\aR^n$ be a Greenian domain (see \cite{AG}), i. e. a domain such that for each $y\in D$ there is the Green function $G_D(x,y)$. Since the Green functions are continuous in both variables, the mapping $\Phi:\,D\to SH^-_1(D,\phi)$ defined as $\Phi(y)=\wtl G_D(x,y)=l^{-1}(y)G_D(x,y)$,
where $l(y)=\|G_D(x,y)\|_{\phi}$, is a homeomorphism on its image. The closure $\wtl D_\phi$ of $\Phi(D)$ in $SH^-_1(D,\phi)$ is compact and consists of $\Phi(D)$ and the set $\bd_MD$ of the limits in $L^1(D,\phi)$ of sequences of functions $\{\wtl G_D(x,y_j)\}$ such that the sequence $\{y_j\}$ has no accumulation points in $D$.   Since these limits are harmonic on $D$, the sets $\bd_MD$ and $\Phi(D)$ do not meet. By Lemma \ref{L:en} if $\phi$ and $\psi$ are norming functions on $D$, then the sets $\wtl D_\phi$ are homeomorphic to each other.
\par In \cite{M} R. S. Martin defines the Martin compactification of $D$ by choosing a point $x_0\in D$ and then adding to $D$ all equivalence classes of converging uniformly on compacta sequences of function $m^{-1}(y_j)G_D(x,y_j)$, where $m(y_j)=G_D(x_0,y_j)$ and the sequence $\{y_j\}$ has no accumulation points. By Harnack's inequality, the third property of a norming function $\phi$ and Lemma \ref{L:ei} the uniform convergence on compacta of harmonic functions is equivalent to convergence in $L^1(D,\phi)$. Therefore, $\wtl D_\phi$ is homeomorphic to the Martin compactification of $D$ and $\bd_MD$ is the Martin boundary of $D$.
\section{The existence of norming functions}
\par Now we will show that every domain $D\sbs\aR^n$ has a norming function. We start with a lemma. Let $B(a,r)$ be the ball of radius $r$ centered at $a$ in $\aR^n$
\bL\label{L:afb} If $D\sbs\aR^n$ is a domain with $C^2$-boundary,
then the function $\phi(x)\equiv1$ on $D$ is norming.\eL
\begin{proof} First we prove that both Green and Poisson kernels are
uniformly integrable on $D$. We may assume that $B(0,s)\sbs D\sbs B(0,1)$ for some $s>0$.
For $\eps>0$ we let $D_\eps=\{x\in D:\,\dl(x)<\eps\}$,
\[v_\eps(y)=\il {D_\eps}G_D(x,y)\,dx\text { and } h_\eps(y)=\il {D_\eps}P_D(x,y)\,dx.\] Since $D\sbs B(0,1)$,
\[G_D(x,y)\ge-\frac{1}{|x-y|^{n-2}}, n\ge 3\text{ and } G_D(x,y)\ge\ln\frac{|x-y|}2, n=2.\]
Thus the function $v_\eps(y)$ is defined and continuous on $D$.
\par Let $n\ge3$. For a point $y\in D_\eps$ take a
ball $B=B(y,r)$ of the radius $r=\dl(y)$. Since
$\dl(x)\le\dl(y)+|x-y|$, by Theorem
\ref{T:kls}
\[G_D(x,y)\ge -A\left(\frac{r^2}{|x-y|^n}+\frac{r}{|x-y|^{n-1}}\right).\]
Therefore
\[v_\eps(y)\ge-\il {B}\frac1{|x-y|^{n-2}}\,dx-
A\il {D_\eps\sm B}\frac{r^2}{|x-y|^n}\,dx-\il {D_\eps\sm
B}\frac{r}{|x-y|^{n-1}}\,dx.\] The first integral is
equal to $c_nr^2$. The second integral does not exceed
\[A\il {B(y,2)\sm B}\frac{r^2}{|x-y|^n}\,dx=c_nr^2\ln\frac2r.\]
The function $|x|^{1-n}$ is integrable on $B(0,2)$. Since the measure of $D_\eps$ converges to 0 as $\eps\to0$, by the absolute continuity of the integral there is $C_1(\eps)\to0$ as $\eps\to0$ such that
\[\il {D_\eps\sm B}\frac{1}{|x-y|^{n-1}}\,dx\le C_1(\eps).\]
Thus
$v_\eps(y)\ge-C(\eps)\dl(y)$ when $y\in D_\eps$, where $C(\eps)\to0$
as $\eps\to0$. Since $v_\eps$ is harmonic on $D'_\eps=D\sm D_\eps$ we
see that
\[v_\eps(y)\ge -C(\eps)\eps\ge-C(\eps)\dl(y)\] on $D$.
In particular, when $\eps=1$
\begin{equation}\label{e:bn4}
v_1(y)=\il DG_D(x,y)\,dx\ge -C(1)\dl(y).
\end{equation}
\par If $y\in\bd D$, then $\dl(x)\le|x-y|$ and by (\ref{e:bn3})
\[h_\eps(y)=\il {D_\eps}P_D(x,y)\,dx\le
\il {D_\eps}\frac{B}{|x-y|^{n-1}}\,dx\le M(\eps),\] where
$M(\eps)$ depends only on $D$ and $\eps$ and converges to $0$ as
$\eps\to0$.
\par The function $G_D(x,y)$, considered as a mapping into the extended real line $[-\infty,\infty]$, is continuous on $D\times D$. Hence there is a constant $d_1(\eps)>0$ such that $G_D(x,y)<-d_1(\eps)$ on $B(0,s)$ when $x,y\in D'_\eps$. Let $h(y)$ be a continuous subharmonic function on $D$ equal to $-d_1(\eps)$ on $B(0,s)$, zero on $\bd D$ and harmonic on $D\sm\ovr B(0,s)$. By the maximum principle, the Keldysch--Lavrentiev--Hopf lemma and (\ref{e:bn4}) we have
\[ G_D(x,y)\le h(y)<-K(\eps)\dl(y)\le d(\eps)v_1(y),\]
where $K(\eps)>0$, $d(\eps)=C(1)K(\eps)$, $x\in D'_\eps$ and $y\in D$.
Also
\[P_D(x,y)=\frac{\bd G_D(x,y)}{\bd n(y)}\ge K(\eps)>0,\]
when $x\in D'_\eps$.
\par Let $u$ be a negative subharmonic function on $D$. By Fubini's theorem
and (\ref{e:pjb})
\beq\label{e:ifb}\il D u\,dy=
\il D v_1\Dl u-\il {\bd D}h_1\,d\mu.\end{equation}
\par The second integral in the right side is finite because
$0\le h_1(y)\le M(1)$. Now for $x\in D'_\eps$ we have
\[\il D v_1\Dl u\ge
\frac1{d(\eps)}\il DG_D(x,y)\Dl u(y)\ge
\frac1{d(\eps)}u(x).\] Therefore, if the first integral in the
right side of (\ref{e:ifb}) is infinite, then $u$ will be equal to
$-\infty$ everywhere. Thus $\|u\|_1<\infty$ and this proves the first property of norming functions.
\par  If $x\in D'_\eps$, then the estimate for $G_D$ above yields
\[\il DG_D(x,y)\Dl u(y)\le d(\eps)\il D v_1\Dl
u.\] The estimate for $P$ tells us that
\[\il{\bd D}P_D(x,y)\,d\mu(y)\ge K(\eps)\il{\bd
D}d\mu\ge\frac{K(\eps)}{M(1)}\il {\bd D}h_1(y)\,d\mu.\] Let
$a(\eps)=\min\{d(\eps),K(\eps)/M(1)\}$. For $x\in D'_\eps$ by (\ref{e:pjb}) and (\ref{e:ifb}) we get
\[u(x)\le a(\eps)\left(\il D v_1(y)\Dl u(y)-
\il{\bd D}h_1(y)d\mu\right)= a(\eps)\il D u(y)\,dy\]
and this proves the second property of norming functions.
\par The function $v_1(y)$ is negative and subharmonic on $D$. Hence by
Keldysch--Lavrentiev--Hopf lemma $v_1(y)<-\al\dl(y)$, $\al>0$. Also
$h_1(y)>\beta>0$ on $\bd D$. Thus
\begin{equation}\begin{aligned} &\il
{D_\eps}u\,dy=\il Dv_\eps\Dl u- \il
{\bd D}h_\eps\,d\mu\ge-C(\eps)\il D\dl\Dl u-
M(\eps)\il{\bd D} d\mu\ge\notag\\&\frac{C(\eps)}\al\il Dv_1\Dl
u-\frac{M(\eps)}\beta\il{\bd D}h_1d\mu\ge b(\eps)\il
Du(y)\,dy,\notag
\end{aligned}\end{equation}
where $b(\eps)=\max\{C(\eps)\al^{-1},M(\eps)\beta^{-1}\}$ converges to 0 as $\eps\to0$. This shows the third property of norming functions.
\par The case $n=2$ has a completely analogous proof.
\end{proof}
\par This lemma fails for general bounded domains as the following
example shows. \bE Let $D=\{z=x+iy\in{\mathbb C}:\,|z|<1, 0<\Arg
z<\pi/2\}$ and let
\[u(z)=\im\frac1{z^2}=-\frac{2xy}{(x^2+y^2)^2}.\]
Then $u$ is a negative harmonic function on $D$ and it is not
integrable.
\eE
\bT\label{T:eaf} Every domain $D\sbs{\mathbb R}^n$ has a norming
function.\eT
\begin{proof} Let $D_j\Sbs D_{j+1}\Sbs D$ be a sequence of
subdomains with smooth boundaries such that $\cup D_j=D$ and let
$G_j=D_{j+1}\sm D_j$.\par
By Lemma \ref{L:afb} there are positive constants $c_j$ such that
\begin{equation}\label{e:cj}
u(x)\le c_j\int\limits_{D_{j+1}}u\,dx
\end{equation}
for each negative subharmonic function $u$ on $D$ and every point $x\in
D_j$. Consequently, there are constants $b_j$, $0<b_j<1$, such that \[\il{D_j}u\,dx\le b_j\il{D_{j+1}}u\,dx\]  and
\[\int\limits_{D_1}u\,dx\le d_j\int\limits_{D_j}u\,dx,\]
where $0<d_j=b_1b_2\dots b_{j-1}<1$.
\par Let $G_j=D_{j+1}\sm D_j$ and let $\phi(x)$ be a positive continuous function on $D$ such that $\phi\le 1$ on $D_1$ and $\phi(x)\le2^{-j}d_j$ on $G_j$.
\par Now \begin{equation}\label{e:ei}\il{D\sm
D_j}u\phi\,dx=
\sum\limits_{k=j}^\infty\il{G_k}u\phi\,dx\ge
\sum\limits_{k=j+1}^\infty d_k2^{-k}\int\limits_{D_k}u\,dx\ge
2^{-j}\int\limits_{D_1}u\,dx.\end{equation}
So
\[\il D u\phi\,dx\ge 2\il {D_1} u\,dx>
-\infty.\] By (\ref{e:ei})
\begin{equation}\label{e:ei2}
\il{D\sm D_j}u\phi\,dx\ge
\frac{2^{-j+1}}c\int\limits_Du\phi\,dx,\end{equation}
 where $c=\inf\phi(x)$ on $D_1$.
\par By (\ref{e:cj}) and (\ref{e:ei2}) for all points $x\in D_j$ we
have
\begin{equation}\label{e:ei3}\begin{aligned}u(x)\le
c_j\il {D_{j+1}}u\,dx\le c_j\il{D_1}u\,dx
\le\frac{c_j}2\il Du\phi\,dx.
\end{aligned}\end{equation} Formulas (\ref{e:ei2}) and
(\ref{e:ei3}) show that the function $\phi$ is
norming.\end{proof}
\section{Norming volume forms on complex manifolds}
\par If $V$ is a positive continuous volume form on a complex manifold $M$, then $L^1(M,V)$ is the space of all Lebesgue functions $u$ on $M$ such that
\[\|u\|_V=\il M|u|\,dV<\infty.\]
\par A positive continuous volume form $V$ on $M$ is {\it norming } if for each compact set $F\Sbs M$ there is a
positive constant $C(F)$ and for every increasing sequence of
open sets $M_j\Sbs M$ with $\cup M_j=M$ there is a
sequence of numbers $\eps_j>0$ converging to zero such that for every
negative plurisubharmonic function $u$ on $M$:
\be\item $\|u\|_V<\infty;$
\item $u(z)\le -C(F)\|u\|_V$ for every point $z\in F$;
\item
$\|u\|_{L^1(M\sm M_j,V)}\le\eps_j\|u\|_V$.
\ee
\bT\label{T:eavf}Every connected complex manifold $M$ has a norming volume form.
\eT
\begin{proof} First, we prove this theorem when $M$ is a relatively compact connected open set with smooth boundary in a complex manifold. Let us take a finite open cover of $\ovr M$ by biholomorphic images $F_j(D'_j)$ of domains $D'_j\sbs\aC^n$, where $n=\dim M$ and $1\le j\le m$. We may assume that the sets $D_j=F_j^{-1}(M\cap F_j(D_j))$ are domains. Then the open sets $U_j=F_j(D_j)\sbs M$ form a finite open cover of $M$. Let $\phi_j$ be a norming function on $D_j$, $G_j=F_j^{-1}$ and $V_j=G_j^*(\phi_j\,dx)$ be the pull-back of the volume form $\phi_j\,dx$ on $D_j$ to $U_j$.
\par Let $\psi_j$ be a partition of unity subordinated to the cover $\{U_j\}$. We let $V=\sum_{j=1}^k\psi_jV_j$. We assume that for all $j$ the sets $\{\psi_j>0\}$ are non-empty. If $u$ is a negative plurisubharmonic function on $M$, then
\[\il Mu\,dV=\sum_{j=1}^m\il {D_j}u(F_j(x))\phi_j(x)\,dx>-\infty.\]
\par Suppose that the intersection of $U_j$ and $U_k$ is non-empty. Let us take a compact set $A\sbs U_j\cap U_k$ such that
\[V_j(A)=\il AdV_j>0.\] There is a constant $c_k>0$ such that
\[u(z)<c_k\il {U_k}u\,dV_k=c_k\il {D_k}u(F_k(x))\phi_k(x)\,dx\]
for any point $z\in A$.
Hence
\[\il {U_j}u\,dV_j\le\il Au\,dV_j<c_kV_j(A)\il {U_k}u\,dV_k.\] This means that there are constants $c_{jk}>0$ such that
\[\il {U_j}u\,dV_j\le c_{jk}\il {U_k}u\,dV_k\] whenever $U_j\cap U_k\ne\emptyset$.
\par Since $M$ is connected for any $1\le j,k\le m$ there is a finite chain of sets $U_{i_1},\dots,U_{i_p}$ such that $U_j=U_{i_1}$, $U_k=U_{i_p}$ and $U_{i_l}\cap U_{i_{l+1}}\ne\emptyset$. Hence for any any $1\le j,k\le m$ there are  constants $c_{jk}>0$ such that
\[\il {U_j}u\,dV_j\le c_{jk}\il {U_k}u\,dV_k.\]
This, in its turn, implies that for any $1\le j\le m$ there is a constant $d_j>0$ such that
\[\il {U_j}u\,dV_j\le \sum_{k=1}^mc_{jk}\il {U_k}\psi_ku\,dV_k
\le d_j\il Mu\,dV.\]
\par Let $F$ be a compact set in $M$. For every point $z\in F\cap U_j$ we take relatively compact open sets $F_z\sbs U_j$ containing $z$ and then choose a finite cover of $F$ by such sets. Let $A_j\Sbs U_{k(j)}$ be the elements of this cover. If $z\in A_j\cap F$, then there is a constant $a_j>0$ such that
\[u(z)<a_j\il {U_{k(j)}}u\,dV_{k(j)}\le a_jd_{k(j)}\il Mu\,dV.\]
Taking as $C(F)$ the minimal constant $a_jc_{k(j)}$ we see that $V$ satisfies the second property of norming volume forms.
\par Let $M_k$ be an increasing sequence of open sets $M_k\Sbs M$ with $\cup M_k=M$. Then for each $j$ there is a sequence of numbers $\eps_{jk}$ converging to zero such that
\[\il{U_j\sm M_k}u\,dV\ge\eps_{jk}\il {U_j}u\,dV_j.\]
Hence
\[\il {M\sm M_k}u\,dV=\sum_{j=1}^m\il {U_j\sm M_k}\psi_ju\,dV_j\ge \sum_{j=1}^m\il {U_j\sm M_k}u\,dV_j\ge \sum_{j=1}^m\eps_{jk}\il {U_j}u\,dV_j.\]
For each $j$ there is a compact set $B_j\sbs U_j$ and a constant $a_j>0$ such that $V_j(B_j)>0$ and $\psi_j>a_j$ on $B_j$.
By Lemma \ref{L:ei} there are constants $b_j>0$ such that
\[\il {U_j}u\,dV_j\ge b_j\il{B_j}u\,dV_j\ge \frac{b_j}{a_j}\il {U_j}\psi_ju\,dV_j.\]
Hence
\[\il {M\sm M_k}u\,dV\ge\sum_{j=1}^m\frac{\eps_{jk}b_j}{d_j}\il {U_j}\psi_ju\,dV_j\ge\eps_k\il Mu\,dV,\]
where $\eps_k$ is some positive sequence converging to 0. This shows the existence of norming forms on relatively compact connected open sets with smooth boundary in a complex manifold.
\par In the general case we exhaust $M$ by relatively compact connected open sets $M_j$ with smooth boundary and repeat the proof of Theorem \ref{T:eaf}.
\end{proof}
\par The first important feature of norming volume forms $V$ on a
connected complex manifold $M$ is the fact that every non-positive plurisubharmonic function
$u\not\equiv-\infty$ is integrable with respect to the measure
$dV$. Hence, the cone $PSH^-(M)$ of all negative plurisubharmonic
functions on $M$  belongs to the Banach space $L^1(M,V)$.
\par Repeating the proof of Lemma \ref{L:en} we get
\bL\label{L:enp} Norming volume forms determine equivalent norms on
the cone of negative plurisubharmonic functions.
\eL
\par Analogously the following plurisubharmonic version of Lemma \ref{L:ei} is valid.
\bL\label{L:eip} Let $V$ be a norming volume on a complex manifold $M$
and let $F$ be a compact set in $M$ with $V(F)>0$. Then there is a positive constant $c$, depending only on $F$, such that for every negative plurisubharmonic function $u$ on $M$ we have
\[c\|u\|_V\le-\il Fu\,dV.\]\eL
\par Let us denote by $B_V$ the closed unit ball in $L^1(M,V)$.
\bT\label{T:ctp} The set $PSH^-_1(M,V)=PSH^-(M)\cap B_V$ is
compact in $L^1(M,V)$.\eT
\begin{proof} Let  $v_j$ be a sequence of plurisubharmonic functions in $PSH^-_1(M,V)$ and let $U\sbs M$ be the biholomorphic image of a domain $D\sbs\aC^n$ by a mapping $F$. Then the functions $u_j=v_j\circ F$ are subharmonic on $D$. By \cite[Theorem 4.1.9]{Ho} there is a subsequence $u_{j_k}$ converging  in $L^1_{loc}(D)$ to a subharmonic function $u$ on $D$. Moreover,
\[w(z)=\limsup_{k\to\infty}u_{j_k}(z)\le u(z)\]
on $D$ and the left and the right side are equal a.e. Thus $u$ is the upper semicontinuous regularization of $w$ and by \cite[Prop. 2.9.17]{K2} $u$ is plurisubharmonic.
\par It follows that if $M_i\Sbs M$ is an increasing sequence of connected open sets such that $\cup M_i=M$, then there is a subsequence $\{v_{j_k}\}$ converging in $L^1(M_i,V)$ for any $i$ to a plurisubharmonic function $v$ on $M$. By the third property of norming volume forms there is a sequence of numbers $\eps_i>0$ converging to zero such that
\[\il{M\sm M_i}u\,dV\ge\eps_i\il Mu\,dV\]
for all $u\in PSH^-(M)$. Hence
\begin{equation}\begin{aligned}
&\il M|v_{j_k}-v|\,dV\le\il{M_i}|v_{j_k}-v|\,dV-\il{M\sm M_i}(v_{j_k}+v)\,dV
\\&\le\il{M_i}|v_{j_k}-v|\,dV-\eps_i\il M(v_{j_k}+v)\,dV\le \il{M_i}|v_{j_k}-v|\,dV+2\eps_i.\notag
\end{aligned}\end{equation}
Thus this subsequence converges to $v$ in $L^1(M,V)$.
\end{proof}
\bP\label{P:tp} If $F:\,M\to N$ is a proper holomorphic mapping between complex manifolds $M$ and $N$ and $V$ is a norming volume form on $N$, then $V^*=F^*V$ is a norming volume form on $M$.
\eP
\begin{proof} Let $A$ be a singular set of $F$ and $A'=F(A)$. Any point in $N\sm A'$ has the same finite number $m$ of preimages under the mapping $F$. If $u$ is a negative plurisubharmonic function on $M$, then $u^*(w)=\sum_{F(z)=w}u(z)$ is a plurisubharmonic function on $N\sm A'$ that is locally bounded above near any point of $A'$. Since the set $A'$ is analytic, $u^*$ extends uniquely to $N$ as a plurisubharmonic function. Thus
\[\il Mu\,dV^*=m\il Nu^*\,dV>-\infty,\]
when $u$ is a negative plurisubharmonic function on $M$.
\par If $G$ is a compact set in $M$, then $G'=F(G)$ is a compact set in $N$ and \[u(z)\le u^*(F(z))\le c\il Nu^*\,dV=\frac cm\il Mu\,dV^*.\]
If $M_j\Sbs M$ is an exhaustion of $M$ and $M'_j=M\sm M_j$, then taking into account that $F$ is proper we get
\[\il {M'_j}u\,dV^*\ge m\il {F(M'_j)}u^*\,dV\ge m\eps_j\il Nu^*\,dV=\eps_j\il Mu\,dV^*,\]
where $\eps_j>0$ is some sequence converging to 0.
\end{proof}
\section{Pluripotential compactification}
\par Let $M$ be a complex manifold. For $w\in M$  we consider the pluricomplex Green function, introduced in \cite{K1},
\[g_M(z,w)=\sup u(z),\]
where the supremum is taken over all negative plurisubharmonic functions $u$ such that the function $u(z)-\log\|z-w\|$ is bounded above near $w$. It is known that $g_M(z,w)$ is plurisubharmonic in $z$. (Here we assume that $u\equiv-\infty$ is a plurisubharmonic function.)
\par The function $g_M(z,w)$ is also maximal in $z$ outside $w$, i.e., if $G\sbs M$ is a domain whose closure does not contain $w$ and $v$ is a plurisubharmonic function on a neighborhood $U$ of $\ovr G$ such that $v(z)\le g_M(z,w)$ on $\bd G$, then $v(z)\le g_M(z,w)$ on $G$. Indeed, if $v$ is a negative plurisubharmonic function on $M$ which is less than $g_M$ on a neighborhood of the boundary of a domain $G\sbs M$, $w\not\in\ovr G$, then we take the function $v_1$ equal to $g_M$ on $M\sm G$ and to $\max\{g_M,v\}$ on $G$. This function will be negative and plurisubharmonic on $M$ and $v_1=g_M$ near $w_0$. Thus $g_M\ge v$ on $G$.
\par We introduce {\it locally uniformly pluri-Greenian} complex manifolds $M$, where every point $w_0\in M$ has a coordinate neighborhood $U$ with the following property:  there is an open set $W\sbs U$ containing $w_0$ and a constant $c$ such that $g_M(z,w)>\log\|z-w\|+c$ on $U$ whenever $w\in W$;
\par If $M$ is a ball $B(w_0,r)$ of radius $r$ centered at $w_0\in\aC^n$, then $g_M(z,w_0)=\log(\|z-w_0\|/r)$. Since $g_M$ is monotonic in $M$, it follows that if $M$ is a bounded domain in $\aC^n$, then $g_M(z,w)\ge\log(\|z-w\|/r)$,
where $r$ is the radius of circumscribed ball of $M$ centered at $w$. Hence bounded domains in $\aC^n$ are locally uniformly pluri-Greenian.
\par We will need a version of \cite[Lemma 6.2.4]{K2}.
\bL\label{L:kl} If $M$ is a locally uniformly pluri-Greenian complex manifold  and $w_0\in M$, then for any $\eps>0$ and any neighborhood $X$ of $w_0$ there is a neighborhood $Y$ of $w_0$ such that
\[1-\eps\le\frac{g_M(z,w_0)}{g_M(z,w)}<1+\eps\]
whenever $w\in Y$ and $z\in M\sm X$.
\eL
\begin{proof} Let $U$ be a coordinate neighborhood of $w_0$ from the definition of locally uniformly pluri-Greenian manifolds. We may assume that $U\sbs X$. By this definition $g_M(z,w)>\log\|z-w\|+c$ on $U$ when $w\in B(w_0,r)\sbs U$ for some $r>0$. On the other hand, if $w\in B(w_0,r/4)$, then $B(w_0,r/2)\sbs B(w,3r/4)\sbs B(w_0,r)$ and by monotonicity of pluricomplex Green functions there is a constant $c_1$ depending only on $r$ such that $g_M(z,w)\le\log\|z-w\|+c_1$ on $B(w_0,r/2)$.
\par If $0<t<r/4$, and $\|w-w_0\|<t/2$, then $\log t+c-2\le g_M(z,w)\le\log t+c_1+2$
on $\bd B(w_0,t)$. Hence there is $0<t_0<r/4$ such that $(1+\eps)g_M(z,w)\le g_M(z,w_0)\le(1-\eps)g_M(z,w_0)$  on $\bd B(w_0,t_0)$ when $w\in B(w_0,t_0/2)$. Our lemma follows with $Y=B(w_0,t_0/2)$ by the maximality of $g_M$.
\end{proof}
\par Let $V$ be a norming volume form on $M$. Let
$c_V(w)=\|g_M(z,w)\|_V$.
We define the mapping $\Phi_V:\,M\to PSH^-_1(M)$ as $\Phi_V(w)=\wtl g_M(z,w)=c^{-1}_V(w)g_M(z,w)$.
\bL If $M$ is a locally uniformly pluri-Greenian complex manifold, then the mapping $\Phi_V$
has the following properties:\be
\item $\Phi_V$ is a continuous bijection onto $\Phi_V(M)$;
\item for every compact set $N\sbs M$ the mapping $\Phi_V$ is a homeomorphism between $N$ and $\Phi_V(N)$.\ee
\eL
\begin{proof} From properties of locally uniformly pluri-Greenian complex manifolds it follows immediately that $\Phi_V$ is a bijection. It follows from Lemma \ref{L:kl} and the inequality $g_M(z,w)>\log\|z-w\|+c$ near $w$ that the function $c_V(w)$ is continuous and, consequently, $\Phi_V$ is continuous.
\par If a set $N\Sbs M$, then $\Phi_V$ is continuous and bijective on $N$. If a sequence $\wtl g_M(z,w_j)$, $w_j\in N$, converges to $\wtl g_M(z,w_0)$ in $L^1(M,V)$, then we take any subsequence $w_{j_k}$ of $\{w_j\}$ converging to $x_0\in N$. By continuity of $\Phi_V$ the sequence $\wtl g_M(z,w_{j_k})$ converges to $\wtl g_M(z,x_0)$ and this implies that $x_0=w_0$. Thus the sequence $\{w_j\}$ converges to $w_0$ in $M$.
\end{proof}
\par The norm of $\Phi_V(w)$ in $L^1(M,V)$ is equal to 1. Hence, by
Theorem \ref{T:ctp} the closure $\wtl M_V$ of $\Phi_V(M)$ in
$L^1(M,V)$ is compact and we call the set $\wtl M_V$ the {pluripotential} compactification of $M$. The set $\wtl M_V$ consists of $\Phi_V(M)$ and the set $\bd_MD$ of the limits in $L^1(M,V)$ of sequences of functions
$\{\wtl G_D(x,y_j)\}$ such that the sequence $\{y_j\}$ has no accumulation points in $D$.  The closure $\bd_P$ of the set $\wtl M_V\sm\Phi_V(M)$ is called the {\it pluripotential} boundary of $M$.
\bL\label{L:uc} Let $M$ and $N$ be locally uniformly pluri-Greenian complex manifolds and $F:\,M\to N$ be a biholomorphism. Let $V$ and $U$ be norming volume forms on
$M$ and $N$ respectively. Then there is a canonical homeomorphism $H$ of $\wtl N_U$ onto $\wtl M_V$ such that $H(\Phi_U(N))=\Phi_V(M)$.\eL
\begin{proof} We define $\wtl H:\,PSH^-_1(N,U)\to PSH^-_1(M,V)$ as $\wtl H(u)=l^{-1}(v)v$, where $v=u\circ F$ and
$l(v)=\|v\|_V$. If we prove that $\wtl H$ is a homeomorphism and $\wtl H(\Phi_U(N))=\Phi_V(M)$, then the restriction $H$ of $\wtl H$ to $\wtl N_U$ will be the required mapping.
\par First of all, we note that $\wtl H$ is bijective. Secondly, if $U^*=F^*U$ and the mapping $\wtl P:\,PSH^-_1(N,U)\to PSH^-_1(M,U^*)$ is defined as $\wtl P(u)=u\circ F$, then $\wtl P$ is a bijective isometry. Finally, by Lemma \ref{L:eip} the function $l(v)$ is continuous on the compact set $PSH^-_1(M,U^*)$. Hence the mapping $v\to l^{-1}(v)v$ is a homeomorphism of $PSH^-_1(M,V^*)$ onto $PSH^-_1(M,V)$. The composition of two latter mappings is $\wtl H$ and our lemma is proved.
\end{proof}
\par In particular, all pluripotential compactifications are homeomorphic to each other and we will denote them by $\wtl M$. Another immediate consequence of this lemma is
\bT Let $M$ and $N$ be locally uniformly pluri-Greenian complex manifolds. Then any biholomorphic mapping $F:\,M\to N$ extends to a homeomorphism of $\wtl M$ onto $\wtl N$. \eT
\section{Examples}
\par When working with examples it is useful to choose a better normalizing factor for pluricomplex Green function. The factor $\|g_M(z,w)\|^{-1}_V$ was optimal for the proofs but hard to calculate in concrete cases. However, if a sequence $\al(w_j)g_M(z,w_j)$ converges in $L^1(M,V)$ to some non-zero function $u$, then the sequence $\wtl g_M(z,w_j)$ also converges to a scalar multiple of $u$.
\par \bE Let $M=\aB^2$ be the unit ball in $\mathbb C^2$.
Evidently, $g_M(z,0)=\log||z||$ and
$g_M(z,a)=g_M(f(z),f(a))=\log||f(z)||$, where $f$ is an automorphism
of the ball transforming $a=(a_1,a_2)$ into 0. If
$w=(w_1,w_2)=f(z)=f(z_1,z_2)$ then
\[w_1=\frac{r-z'_1}{1-z'_1r};\qquad
w_2=\frac{\sqrt{1-r^2}}{1-z'_1r}z'_2,\]where
$r=||a||,\,z'_1=r^{-1}(z,a)$ and $z'_2=r^{-1}(z,\bar a)$.
Therefore,
\[g_M(z,a)=\log\left|1-\frac{(1-||z||^2)(1-r^2)}{|1-z'_1r|^2}\right|.\]
As a normalizing factor we take $|g(0,a)|=-2\log r$. Then
\[\wtl g_M(z,a)=-\frac{1}{2\log r}\log\left|1-
\frac{(1-||z||^2)(1-r^2)}{|1-z'_1r|^2}\right|.\]
If a sequence $\wtl g_M(z,a_j)$ converges and $\|a_j\|\to1$, then $a_j\to a=(a_1,a_2)\in \bd\aB^2$ and the limit is
\[\wtl g_M(z,a)=\frac{||z||^2-1}{|1-(z,a)|^2}.\] The
function $\wtl g_M(z,a)$ is maximal because for mappings
$f(\zeta)=(\zeta,C(1-\zeta)):\aD\to \aB^2$ the functions $\wtl g_M(f(\zeta),a)$ are harmonic.
\par So in this case $\wtl M=\ovr \aB^2$,  $\bd_PM=\bd\aB^2$ and the mapping $\Phi_V$ is a homeomorphism of $\ovr \aB^2$ onto $\wtl M$. This means that the Euclidean boundary and the pluripotential boundary coincide.
\eE
\par \bE  Let $M=\aD^2\sbs\aC^2$, $z=(z_1,z_2)$ and $w=(w_1,w_2)$. Then
\[g_M(z,w)=\log\max\left\{\left|\frac{z_1-w_1}{1-z_1\bar w_1}\right|,\left|\frac{z_2-w_2}{1-z_2\bar w_2}\right|\right\}.\]
As a normalizing factor we take $\al(w)=|g^{-1}_M(0,w)|=-\log^{-1}\max\{|w_1|,|w_2|\}$.
\par If a sequence $w_j=(w_{1j},w_{2j})$ in $M$ converges to $w_0=(w_{10},w_{20})$ and $|w_{10}|=1$ while $|w_{20}|<1$, then the sequence $\wtl g_M(z,w_j)$ converges to
\begin{equation}\label{E:bd1}\wtl g_M(z,w_0)=\frac{|z_1|^2-1}{|1-z_1\bar w_{10})|^2}.\end{equation}
Similarly, if $|w_{01}|<1$ while $|w_{02}|=1$, then the sequence $\wtl g_M(z,w_j)$ converges to
\begin{equation}\label{E:bd2}\wtl g_M(z,w_0)=\frac{|z_2|^2-1}{|1-z_2\bar w_{20})|^2}.\end{equation}
\par If $|w_{01}|=|w_{02}|$, then the sequence $\wtl g_M(z,w_j)$ converges if and only if the sequence $(\log|w_{1j}|)/\log|w_{2j}|$ has the finite or infinite limit $c$. If $c=\infty$, then $\wtl g_M(z,w_j)$ converges to the function from (\ref{E:bd2}), while if $c=0$, then the limit of $\wtl g_M(z,w_j)$ is the function from (\ref{E:bd1}).
\par If $0<c\le1$, then $\wtl g_M(z,w_j)$ converges to the function
\[\wtl g_M(z,w_0)=\log\max\left\{\frac{|z_1|^2-1}{|1-z_1\bar w_{10})|^2},c^{-1}\frac{|z_2|^2-1}{|1-z_2\bar w_{20})|^2}\right\}.\]
If $1\le c<\infty$, then $\wtl g_M(z,w_j)$ converges to the function
\[\wtl g_M(z,w_0)=\log\max\left\{c\frac{|z_1|^2-1}{|1-z_1\bar w_{10})|^2},\frac{|z_2|^2-1}{|1-z_2\bar w_{20})|^2}\right\}.\]
\par All limit functions are maximal. The non-distinguished part of $\bd M$ is squeezed into two circles while every point in the distinguished boundary, that is a 2-torus $\aT^2$, blows up to an interval $(0,\infty)$ connecting these circles. If we add to every point of $\aT^2$ the interval $(0,1]$ and the circle from (\ref{E:bd1}) we will get a filled torus in $\aR^3$. Adding interval $(1,\infty)$ and the circle from (\ref{E:bd2}) we will get another filled torus in $\aR^3$. Thus $\bd_PM$ is the double of a filled torus or  the product of a circle and a 2-sphere.
\eE
\par \bE Let $M$ be a smooth strongly convex domain in $\aC^n$. A {\it complex geodesic } is a holomorphic map $\phi:\,\aD\to M$ which is an isometry between the Poincar\'e metric on $\aD$ and the Kobayashi distance $k_M$ on $M$. According to Lempert (see \cite{L1}) on smooth strongly convex domains complex geodesics are injective maps smooth up to the boundary and the Kobayashi and Carath\'eodory distances coincide. The latter implies (see \cite{PS} or \cite{BPT}) that $g_M(z,w)=\log(\tanh k_M(z,w))$.
\par In \cite[Theorem 3]{CHL} the authors constructed a continuous mapping $\Phi:\,\bd M\times M\to\aB^n$ that is smooth on $\bd M\times\ovr M$ outside of the diagonal in $\bd M\times\bd M$ and has the following properties:\be
\item for every $p\in\bd M$ there is $q\in\bd\aB^n$ such that $\Phi(p,\phi(\zeta))$, where $\phi$  is any complex geodesic in $M$ and $p\in\phi(\aT)$, is a complex geodesic in $\aB^n$ passing through $q$;
\item for a fixed $p\in\bd M$ the mapping $\Phi(p,z)$ is a homeomorphism of $\ovr M$ onto $\aB^n$ smooth outside of $\{p\}$.
\ee
\par For a complex geodesic $\phi$ such that $\phi(0)=w$ and $\phi(\zeta)=z$ we choose a point $p\in\phi(\aT)$. Due to the isometry properties of $\Phi$ the value of the function $\log(\tanh k_{\aB^n}(\Phi(p,w),\Phi(p,z))$ does not depend on the choice of $p$ and is equal to $\log(\tanh k_M(z,w))=g_M(z,w)$.
\par Let $\{w_j\}\sbs B$ be a sequence converging to $p\in\bd M$. For a complex geodesic passing through $w_j$ and $z$  we choose as $p_j(z)$ the nearest point to $p$ in $\phi(\aT)$. By the continuity of $\Phi$ the mappings $\Phi(p_j(z),z)$ converge to $\Phi(p,z)$ uniformly on compacta in $M$. Let $x_j\in M$ be the points such that $\Phi(p_j(x_j),x_j)=0$ and let $\lm_j=-2\log^{-1}\|\Phi(p_j(x_j),w_j)\|$. Then the functions
\[\lm_jg_M(z,w_j)=\lm_j\log(\tanh k_{\aB^n}(\Phi(p_j(z),z),\Phi(p_j(z),w_j))\]
converge uniformly on compacta to $\wtl g_M(z,p)=\wtl g_{\aB^n}(\Phi(p,z),p)$.
This function is maximal because it is harmonic on geodesics passing through $p$, equal to 0 on $\bd M\sm\{p\}$ and smooth. So by \cite[Theorem 7.3]{BPT} $\wtl g_M(z,p)$ is a scalar multiple of the function $\Om_{M,p}(z)$. In \cite{BPT} the latter function is called {\it pluricomplex Poisson kernel} of $M$ and it is equal to the derivative of $g_M(z,p)$ along the outside normal at $p$ like in the classical formula (\ref{e:dgf}).
\par In \cite{D} Demailly introduced the notion of pluriharmonic Poisson kernels that depend on the choice of a measure on the boundary and are scalar multiples of each other. It was proved in \cite{BPT} that $\Om_{M,p}$ is a pluriharmonic Poisson kernel in the sense of Demailly. He also computed these kernels for the ball and the polydisk and they are scalar multiples of the functions computed in Examples 1 and 2.
\eE


\begin{thebibliography}{999}
\bibitem{AG}D. H. Armitage, S. J. Gardiner, {\em Classical Potential Theory,} Springer, 2001
\bibitem{BPT} F. Bracci, G. Patrizio, S. Trapani, {\em The pluricomplex Poisson kernel for strongly convex domains,} Trans. Amer. Math. Soc. {\bf 361} (2009),  979–-1005
\bibitem{CHL} C.-H. Chang, M. C. Hu, H.-P. Lee, {\em Extremal analytic discs with prescribed boundary data,} Trans. Amer. Math. Soc. {\bf 310} (1988), 355–-369
\bibitem{D} J.-P. Demailly, {\em Mesures de Monge-Amp\'ere et mesures pluriharmoniques,} Math. Z. {\bf 194} (1987), 519–-564
\bibitem{H} L. L. Helms, {\em Introduction to Potential Theory,}
Wiley-Interscience, 1969
\bibitem{Ho} L.  H\"ormander, {\em The Analysis of Linear  Partial
Differential Operators, I, Distribution Theory and  Fourier
Analysis, } Springer, 1983
\bibitem{K1} M. Klimek, {\em Extremal plurisubharmonic functions
 and invariant pseudodistances,} Bull. Soc. Math. France, {\bf 113} (1985),
123--142
\bibitem{K2} M. Klimek, {\em Pluripotential Theory,} Clarendon Press, 1991
\bibitem{KL} M. V. Keldysch, M. A. Lavrentiev, {\em Sur une
\'evalution pour la fonction de Green,} Dokl. Acad. Nauk USSR, {\bf 24} (1939),
102--103
\bibitem{L} P. Lelong, {\em Discontinuit\'e et annulation de l'op\'erateur de Monge-Amp\'ere complexe,} P. Lelong-P. Dolbeault-H. Skoda analysis seminar, 1981/1983, 219–-224, Lecture Notes in Math., {\bf 1028}, Springer, Berlin, 1983.
\bibitem{L1} L. Lempert, {\em La m\'etrique de Kobayashi et la repr\'esentation des domaines sur la boule,} Bull. Soc. Math. France {\bf 109} (1981), 427–-474
\bibitem{M} R. S. Martin, {\em Minimal positive harmonic functions,}  Trans. Amer. Math. Soc. {\bf 49} (1941), 137–-172
\bibitem{PS} E. A. Poletsky, B. V. Shabat, {\em Invariant metrics,} Current problems in mathematics. Fundamental directions, {\bf 9} (1986),  73–-125,  Itogi Nauki i Tekhniki, Akad. Nauk SSSR
\bibitem{PK} I. I. Privalov, P. K. Kuznetsov, {\em Boundary problems and
classes of harmonic and subharmonic functions in arbitrary domains,}
Mat. Sb., {\bf 6} (1939), 345--375
\bibitem{S} E. D. Solomentsev, {\em Boundary values of subharmonic
functions,} Czech. Math. J., {\bf 8} (1958), 520--534
\bibitem{V} V. S. Vladimirov, {\em Equations of Mathematical Physics,}
Nauka, Moscow, 1967
\bibitem{W} K.-O. Widman, {\em Inequalities for the Green function and boundary continuity of the gradient of solutions of elliptic di®erential equations,} Math. Scand. {\bf 21} (1967), 17--37.
Corresponding lower bound estimates are given in
\bibitem{Z} Z. X. Zhao, {\em Green function for Schrödinger operator and conditioned Feynman-Kac gauge,} J. Math. Anal. Appl. {\bf 116} (1986), 309-–334.

\end{thebibliography}
\end{document}